\documentclass[12pt]{amsart}

\usepackage{amssymb}
\usepackage{epsf}
\usepackage{amsmath,epsf}
\usepackage[latin1]{inputenc}
\usepackage{amsfonts}

\numberwithin{equation}{section}
\newtheorem{thm}{Theorem}[section]

\newtheorem{lem}[thm]{Lemma}
\newtheorem{cor}[thm]{Corollary}
\newtheorem{definition}[thm]{Definition}

\newtheorem{rem}[thm]{Remark}

\def\Q{{\mathbb Q}}

\def\J{{\mathcal J}}

\def\Young#1{\vbox{\smallskip\offinterlineskip
    \halign{&\vbox{##}\kern-\Thickness\cr #1}}}

\newdimen\Squaresize \Squaresize=20pt
\newdimen\Thickness \Thickness=.1pt
\newdimen\Correction \Correction=7pt

\def\Vide#1{\hbox{
       \vbox to \Squaresize{\vss
          \hbox to \Squaresize{\hss#1 \hss}\vss}
    \hskip-\Correction}
   \kern-\Thickness}

\def\Carre#1{\hbox{\vrule width \Thickness
   \vbox to \Squaresize{\hrule height \Thickness\vss
      \hbox to \Squaresize{\hss$\scriptstyle#1$\hss}
   \vss\hrule height\Thickness}
   \unskip\vrule width \Thickness}
   \kern-\Thickness}

\def\Dyck#1{\vbox{\smallskip\offinterlineskip
    \halign{&\vbox{##}\kern-\Thickness\cr #1}}}

\newdimen\Squaresize \Squaresize=15pt
\newdimen\MyDim \MyDim=15pt
\newdimen\Thickness \Thickness=.1pt
\newdimen\Correction \Correction=7pt

\def\Gauche#1{\hbox{\vrule width \Thickness
       \vbox to \Squaresize{\vss
          \hbox to \Squaresize{\hss#1\hss}
       \vss}
    \unskip\kern\Thickness}
   \kern-\Thickness}

\def\Haut#1{\hbox{\kern-\Thickness
       \vbox to \Squaresize{\hrule height \Thickness\vss
          \hbox to \MyDim{\hss#1\hss}
       \vss}
    \unskip}
   \kern-\Thickness}

\def\Bas#1{ \hbox{\kern-\Thickness
       \vbox to \Squaresize{\vss
          \hbox to \MyDim{\hss#1\hss}
       \vss\hrule height\Thickness}
    \unskip}
   \kern-\Thickness}
\def\G{\Gauche{}}

\title[Catalan and Super-Harmonics]{Catalan paths, Quasi-symmetric functions\\ and
Super-Harmonic Spaces}

\author[J.-C.~Aval and N.~Bergeron]
{J.-C.~Aval and N.~Bergeron}

\address[Jean-Christophe Aval]
{Laboratoire A2X\\ Universit\'e Bordeaux 1\\ 351 cours de la
Lib\'eration\\ 33405 Talence cedex\\ FRANCE}

\address[Nantel Bergeron]
{Department of Mathematics and Statistics\\ York University\\
   To\-ron\-to, Ontario M3J 1P3\\ CANADA}

\email[Jean-Christophe Aval]{aval@math.u-bordeaux.fr}
\email[Nantel Bergeron]{bergeron@mathstat.yorku.ca}
\urladdr[Nantel Bergeron]{http://www.math.yorku.ca/bergeron}

\date{\today}
\thanks{N. Bergeron is supported in part by NSERC, PREA and CRC}

\begin{document}

\begin{abstract}
We investigate the 
quotient ring $R$ of the ring of formal power
series
$\Q[[x_1,x_2,\ldots]]$ over the closure of the ideal generated by non-constant
quasi-\break symmetric functions.  We show that a Hilbert basis of the quotient is
naturally indexed by Catalan paths (infinite Dyck paths). We also give a filtration of
ideals related to Catalan paths from
$(0,0)$ and above the line
$y=x-k$. We investigate as well the quotient ring $R_n$ of
polynomial ring in $n$ variables over the ideal generated by non-constant
quasi-symmetric polynomials. We show that the dimension of $R_n$ is bounded above by
the $n$th Catalan number. 
\end{abstract}

\maketitle

%%%%%%%%%%%%%%%%%%%%%%%%%%%%%%%%%%%%%%%%%%%%%%%%
\section{Introduction}

 The ring $Qsym$ of quasi-symmetric functions was introduced by
Gessel~\cite{Ges} as a source
of generating functions for $P$-partitions~\cite{Stanley_EC1}.
Since then, quasi-symmetric
functions have appeared in many combinatorial
contexts~\cite{BMSW,Stanley_EC1,Stanley_ECII}. The relation of $Qsym$ to the ring of
symmetric functions was first clarified by Malvenuto and Reutenauer \cite{MR} via
a graded Hopf duality to the  Solomon descent algebras, then Gelfand {\it
et.~al.~}~\cite{GKal} defined the graded Hopf algebra $NC$ of non-commutative symmetric
functions and identified it with the Solomon descent algebra. In  recent literature, we
see a growing interest in quasi-symmetric functions and non-commutative symmetric
functions as refinements of the ring of symmetric functions.

One unexplored avenue is as an analogue of the (symmetric) harmonic spaces.
A classic combined result of Artin and Steinberg~\cite{Artin, Steinberg} shows that 
the quotient ring of the polynomial ring
$\Q[x_1,x_2,\ldots,x_n]$ in
$n$ variables over the ideal ${\mathcal I}_n$ generated by non-constant 
symmetric polynomials has dimension
$n!$. In fact, this space is a graded symmetric group module that affords the left
regular representation. Moreover under the scalar product 
 $$\langle P, Q \rangle= \big(P(\partial_{x_1},\partial_{x_2},\ldots,\partial_{x_n})
   Q\big)(0,0,\ldots, 0)$$
the graded orthogonal complement $H_n={\mathcal I}_n^\perp$ is a set of representatives 
for the quotient spanned by all possible partial derivatives of the Vandermonde
determinant. We see $H_n$ as the set of polynomials in $n$ variables that are killed by
all symmetric partial derivative operators. In particular, the Laplacian $\sum
\partial_{x_i}^2$ kills any such polynomial, thus the $H-n$ if often called the space of
harmonics. Refinement and generalization of this result has lead to an explosion
of incredible results and conjectures, see
\cite{aval,aval2,ABB,ANB,Berg et
al,Berg2,nantel_adriano,CP,orbit,GH,GP,haiman1,haiman,Tani} for a small portion of this. 
The so-called
$n!$-conjecture (Theorem) of Garsia and Haiman~\cite{GH}, just recently proven by
Haiman~\cite{haiman}, and its connection with Macdonald polynomials~\cite{haiman,Mac} 
is a great achievement in this context.

Here we are interested in the  quotient ring $R$ of the ring of formal power series
$\Q[[x_1,x_2,x_3,\ldots]]$ over the closure of the homogeneous ideal
$\J$ generated by all non-constant quasi-symmetric functions. That is the quotient
  \begin{equation}\label{SPACE}
    R = \Q[[x_1,x_2,x_3,\ldots]] \Big/ \lower2pt\hbox{$ \overline{\J}$}.
   \end{equation}
This quotient is in fact a Hopf algebra. It will be interesting to study its structure
in more detail in future work. Here we concentrate our attention on its linear
structure only.

To every monomial $x_1^{\tilde\alpha_1}x_2^{\tilde\alpha_2}x_3^{\tilde\alpha_3}\cdots$ 
(of finite total degree) in
$\Q[[x_1,x_2,x_3,\ldots]]$ we associate a path in the plane as follow:
  $$ (0,0)\to(\tilde\alpha_1,0)\to(\tilde\alpha_1,1)\to(\tilde\alpha_1+\tilde\alpha_2,1)
     \to(\tilde\alpha_1+\tilde\alpha_2,2)\to(\tilde\alpha_1+\tilde\alpha_2+\tilde\alpha_3,2)
\to\cdots
  $$
If this path remains above the line $y=x$ we say that the path is a {\sl Catalan path}
(or infinite Dyck path). Our main result is 

\begin{thm} \label{main} A monomial Hilbert basis of $R$ is given by the monomials of\break
$\Q[[x_1,x_2,x_3,\ldots]]$ corresponding to Catalan paths.
\end{thm}

We also consider a special filtration of ideals
$\J^{(e)}$ and their respective quotients, such that $\J=\J^{(0)}$ and 
$\J^{(e)}\subseteq\J^{(e+1)}$. The Hilbert basis of each quotient  is indexed by paths
above the line
$y=x-e$. 

This result relates to the harmonic polynomials in the following way. Consider  the
quotient of  the polynomial ring
$\Q[x_1,x_2,\ldots,x_n]$ over the ideal
$\J_n$ generated by all non-constant quasi-symmetric polynomials. Since the ring  of
symmetric polynomials is a subring of the ring of quasi-symmetric polynomials, we have
that ${\mathcal I}_n\subseteq \J_n$. We thus consider the quotient
  \begin{equation}\label{nSPACE}
    R_n = \Q[x_1,x_2,\ldots,x_n] \Big/ \J_n.
   \end{equation}
The space $SH_n=\J_n^\perp$ of representatives for the quotient is a subspace of
harmonic polynomials as
  $$SH_n\subseteq H_n.$$
For this reason we call $SH_n$ the space of {\sl super-harmonic} polynomials. Recall that
$C_n={1\over n+1}{2n \choose n}$ are the famous Catalan numbers. The passage from 
infinitely many variables to finitely many variables is a priori non-trivial, and
requires more work. Here we show:

\begin{thm} \label{Finite}
  \begin{equation}\label{dSPACE}
    \dim SH_n = \dim R_n \le C_n.
   \end{equation}
\end{thm}

In fact we expect equality to hold in (\ref{dSPACE}) and we are still in the
process of completing a proof of this together with Fran\c{c}ois Bergeron and Adriano
Garsia, who have discovered the spaces $SH_n$ and $R_n$ completely independently and in
the same period that we did.  They have conjectured many of the results presented here,
and much more. The results for $SH_n$ and the relations between $H_n$ and $SH_n$ are
extremely interesting, and are the object of an ongoing
collaboration with F.~Bergeron and A. Garsia. We plan to write at least two papers:
one
\cite{deux}  dedicated to a proof of equality in (\ref{dSPACE}), and another
\cite{JCfrananad} to investigate further properties of $SH_n$ and its generalization. 
In particular, the finite version of the successive quotients by the ideals $\J_n^{(e)}$
is related to the work of \cite{Berg2}.  Much of these results can  be
explained in a more general framework and will be the object of further study.  We are
convinced that these results are but the tip of a new iceberg.
 In particular, we would like to find any  natural algebras acting on these
spaces. What are the possible generalizations and specializations of the
super-harmonics?

We underline here that F. Hivert~\cite{hivert} has developed an action of  the Hecke
algebra for which a polynomial is invariant if and only if it is quasi-symmetric.
Unfortunately Hivert's action is not compatible with multiplication and does not
preserve the ideal $\J_n$, hence it does not induce the desired action on the quotient.
It is still interesting to note that Hivert's action is also related to Catalan
numbers.  One way to reformulate his result in~\cite{hivert} is as follows. Consider the
generator $e_i={q-T_i\over(1+q)}$ of the Hecke algebra, where $T_i$ are the standard
generators. Then 
  \begin{equation}\label{temp}
    e_{i}e_{i\pm 1}e_{i} - {q\over (1+q)^2} e_{i}
  \end{equation}
acts, via Hivert's action, as zero on the polynomial ring. Hence, the  quotient of the
Hecke algebra by  Relation~(\ref{temp}), classically known as the Temperley-Lieb
algebra~\cite{jones}, naturally acts on polynomials. This algebra is known to have
dimension equal to $C_n$.

In Section 2 we recall appropriate definitions. In Section 3 we introduce  a special
family of generators for the ideal $\J$ and the associated filtration $\J^{(e)}$. In 
Section~4 we use these generators to show that the monomials corresponding to Catalan
paths span our quotient, as well as the analoguous result for
$\J^{(e)}$. Theorem~\ref{Finite} follows from this section. To complete the proof
Theorem~\ref{main} we use a Gr\"obner basis argument  in Section~5 to show independence.

%%%%%%%%%%%%%%%%%%%%%%%%%%%%
\section{Basic definitions}

A composition $\alpha =[\alpha_1,\alpha_2, \ldots ,\alpha_k]$ of a
positive integer $d$ is an ordered list of positive integers whose sum is
$d$. We denote this by $\alpha \models d$. We   call the integers 
$\alpha_i$ the {\sl parts} of $\alpha$, and denote the number of
parts in $\alpha$ by $\ell(\alpha)$. Given two compositions  $\alpha
=[\alpha_1,\alpha_2, \ldots ,\alpha_k]$ and $\beta =[\beta_1,\beta_2, \ldots
,\beta_\ell]$, we denote by $\alpha\beta$ the concatenation product
$[\alpha_1,\alpha_2, \ldots ,\alpha_k,\beta_1,\beta_2, \ldots ,\beta_\ell]$. Also,
there exists a natural one-to-one correspondence between compositions of $d$ and
subsets of
$\{1,2,\ldots,d-1\}$. If
$A=\{ a_1,a_2,\ldots ,a_{k-1}\}\subset[d-1]$, where  $a_1<a_2<\ldots <a_{k-1}$, then
$A$ corresponds to the composition,
$\alpha =[a_1-a_0,a_2-a_1,\ldots ,a_k-a_{k- 1}]$,  where $a_0=0$ and $a_k=d$. For ease
of notation, we shall denote the set corresponding to a given composition $\alpha$ by
$D(\alpha)$. For compositions $\alpha$ and $\beta$  we say that $\alpha$ is a
\emph{refinement} of
$\beta$  if $D(\beta)\subset D(\alpha)$, and denote this by
$\alpha\preccurlyeq \beta$.

For any composition $\alpha =[\alpha_1,\alpha_2,\ldots ,\alpha_k]$  of $d$  we denote
by $M_\alpha$ the {\sl monomial quasi-symmetric function} \cite{Ges}
  $$M_\alpha(x_1,x_2,\ldots)=\sum_{i_1<i_2<\cdots <i_k} x^{\alpha_1}_{i_1}\ldots 
    x^{\alpha_k}_{i_k}.$$
This is a homogeneous infinite series of degree $d$. We define $M_0=1$,
where $0$ denotes the unique empty composition of $0$.  It is known from the work of
Gessel that the monomial quasi-symmetric functions form a linear basis of a ring (in
fact a Hopf algebra) $Qsym$ of quasi-symmetric functions. 

 An other useful basis of the ring $Qsym_n$ is given by 
the {\sl fundamental quasi-symmetric function}  \cite{Ges}:
  $$F_\alpha(x_1,x_2,\ldots)=\sum_{\alpha\succcurlyeq \beta} M_\beta(x_1,x_2,\ldots) = 
    \sum_{\stackrel{\mbox{\scriptsize $j_1\leq j_2\leq\cdots\leq j_d$}}
                    {i\in D(\alpha)\,\, \Rightarrow\,\, j_i<j_{i+1}}}
    x_{j_1} x_{j_2}\cdots x_{j_d}\,.$$
Fundamental quasi-symmetric functions satisfy the following obvious, but crucial,
relations. For $\alpha =[\alpha_1,\alpha_2,\ldots ,\alpha_k]\models d$, if
$\alpha_1>1$, then
  \begin{equation}\label{firstFrel}
    F_\alpha(x_1,x_2,\ldots) = x_1 F_{[\alpha_1-1, \alpha_2,\ldots,\alpha_k]}(x_1,x_2,\ldots) 
+ F_\alpha(x_2,x_3,\ldots),
   \end{equation}
and if $\alpha_1=1$, then
  \begin{equation}\label{secondFrel}
    F_\alpha(x_1,x_2,\ldots) = x_1 F_{[\alpha_2,\alpha_3,\ldots,\alpha_{k}]}(x_2,x_3,\ldots) 
+ F_\alpha(x_2,x_3,\ldots).
   \end{equation}
Here $F_{\alpha}(x_2,x_3,\ldots)$ is the function $F_{\alpha}(x_1,x_2,\ldots)$
in which the variable $x_i$ is replaced  by $x_{i+1}$. We will see that these relations
are the key ingredients in our proof.

In the following we have to consider  generalized (infinite) compositions.  That is a
sequence
$\tilde\alpha=(\tilde\alpha_1,\tilde\alpha_2,\ldots)$ such that the parts
 $\tilde\alpha_j\ge 0$ for $j\ge 1$ (we allow some
parts to be zero) and the sum of the  parts $d(\tilde\alpha)=\sum
\tilde\alpha_i<\infty$. We say that
$\tilde\alpha$ is a generalized composition of $d(\tilde\alpha)<\infty$.  We use a ``\
$\tilde{\ }$\ '' to indicate that we have a generalized composition and no
``\ $\tilde{\ }$\ '' if the composition  is {\sl standard}, that is without zeros.  We
also consider generalized compositions of finite length and denote by
$\ell(\tilde\alpha)$ the number of parts  of
$\tilde\alpha$. The concatenation of a finite length generalized composition
$\tilde\alpha$ with an infinite one $\tilde\beta$ is denoted by
$\tilde\alpha\tilde\beta$. We also write $\tilde\alpha+\tilde\beta$ and
$\tilde\alpha\le\tilde\beta$ to denote the componentwise sum and componentwise
inequalities, respectively. For an infinite generalized composition  $\tilde\alpha$,
since
$d(\tilde\alpha)<\infty$, only finitely many parts of $\tilde\alpha$ are non-zero. Thus
there is always a finite generalized composition $\tilde\nu$ such that
$\tilde\alpha=\tilde\nu\,0\,0\cdots$.

In this paper, we devote our attention to the ideal $\J=\langle F_\alpha(x_1,x_2,\ldots)
\rangle_{\alpha\models d>0}$ of\break
$\Q[[x_1,x_2,\ldots]]$ generated by the  non-constant quasi-symmetric functions, and
consider  the quotient
  \begin{equation}\label{SPACE}
    R = \Q[[x_1,x_2,\ldots]] \Big/ \lower2pt\hbox{$\overline{\J}$},
   \end{equation}
where $\overline{\J}$ denotes the closure (with respect to the standard topology with
formal power series) of
$\J$ in
$\Q[[x_1,x_2,\ldots]]$.

%%%%%%%%%%%%%%%%%%%%%%%%%%%%
\section{The generators $G_{\tilde\alpha}$} 

In the previous section we noted the relations (\ref{firstFrel}) and (\ref{secondFrel}). 
{}From the first one we deduce that if $\alpha_1>1$, then
  $$F_\alpha(x_2,x_3,\ldots) = F_\alpha(x_1,x_2,\ldots) - x_1
   F_{[\alpha_1-1,\alpha_2,\ldots,\alpha_k]}(x_1,x_2,\ldots).$$
 Since both $F_\alpha(x_1,x_2,\ldots)$ and
$F_{[\alpha_1-1,\alpha_2,\ldots,\alpha_k]}(x_1,x_2,\ldots)$ are in $\J$, we conclude
that
$F_\alpha(x_2,x_3,\ldots)\in\J$. We want to  exploit these properties to a maximum. For
this we construct a set $\{G_{\tilde\alpha}\}\subseteq \J$  indexed by the generalized
(infinite) composition
$\tilde\alpha$ such that 
there exists a factorization $\tilde\alpha=\tilde\pi\tilde\rho$ where
  \begin{equation}\label{condition}
       d(\tilde\pi)-\ell(\tilde\pi)\ge 0.
  \end{equation}

For this, we first define recursively the  functions $G_{\tilde\alpha}$ for all
infinite generalized composition
$\tilde\alpha$. Then in Lemma~\ref{InJ} we
characterize the $\tilde\alpha$ obtained from the transitive closure of the
$G_{\tilde\alpha}\in\J$. Let $\tilde\alpha=\tilde\nu\,0\,0\cdots$
where
$\ell(\tilde\nu)<\infty$ and the last part of
$\tilde\nu$ is non-zero, or $\tilde\alpha=0\,0\cdots$. 
Our definition is recursive on $n=\ell(\tilde\nu)$. If
$\tilde\nu=\nu$ is a standard composition, then let
  \begin{equation}\label{baseG}
    G_{\tilde\alpha}= F_\nu(x_1,x_2,\ldots).
   \end{equation}
If $\ell(\tilde\nu)=0$, then this formula gives $G_{0\,0\,\cdots}=1$.
Assume now that $\tilde\nu$ is non-standard  and let $\tilde\nu=\tilde\gamma\, 0\,
a\,\beta$ be the unique factorization of $\tilde\nu$ such that $a>0$ is a positive
integer,
$\beta$ is a (possibly empty) standard  composition and $\tilde\gamma$ is a (possibly
empty) generalized composition. For $\tilde\alpha=\tilde\gamma\, 0\,
a\,\beta\,0\,0\cdots$ and 
$k=\ell(\tilde\gamma\,0)=\ell(\tilde\gamma)+1$, we define
  \begin{equation}\label{inductG}
    G_{\tilde\alpha}= G_{\tilde\gamma\,  a\,\beta\,0\cdots} - x_k G_{\tilde\gamma\, 
      (a-1)\,\beta\,0\cdots}.
   \end{equation}
Both term on the right are well defined by  induction since
$\ell(\tilde\gamma\,a\,\beta)=\ell(\tilde\nu)-1<n$. 
 
     We now characterize the transitive closure of the definition (\ref{baseG}) and
     (\ref{inductG}) within $\J$. At this  point it is useful to introduce the
following family of ideals. For any $e\ge 0$, let
  $$\J^{(e)} \ =\ \langle F_{\alpha}  \,:\, \exists\, \pi\rho=\alpha, \
d(\pi)-\ell(\pi)\ge e\rangle.
  $$
This is a filtration $\J^{(e)}\subseteq \J^{(e+1)}$ such that $\J=\J^{(0)}$. For a
generalized composition $\tilde\alpha$, we say that it {\sl reaches level $e$} if there
exists a factorization
$\tilde\alpha=\tilde\pi\tilde\rho$ such that
  \begin{equation}\label{econdition}
    d(\tilde\pi)-\ell(\tilde\pi)\ge e.
   \end{equation}
  \begin{lem} \label{InJ} \ 
   \begin{enumerate}
     \item If $\tilde\alpha$ reaches level $e$, then
           $G_{\tilde\alpha}\in\J^{(e)}$. 
     \item Conversely, in (\ref{inductG}), if ${\tilde\gamma\,(a-1)\,
\beta\,0\,0\cdots}$ reaches level $e$, then
$\tilde\alpha$
        reaches level $e$.
   \end{enumerate}
\end{lem}

\proof
For the first statement we proceed by induction on $\ell(\tilde\nu)$ where
$\tilde\alpha=\tilde\nu\,0\,0\cdots$  and the last part of
$\tilde\nu$ is non-zero. If $\ell(\tilde\nu)=0$,  then $G_{0\,0\cdots}=1$ is not in any
of the ideals
$\J^{(e)}$. Assume that $\ell(\tilde\nu)>0$.

We first consider the case when 
$\tilde\nu=\nu$ is a standard composition.  If $\tilde\alpha$ reaches level $e$, then so
is $\nu$ and we have $G_{\tilde\alpha}=F_\nu(x_1,x_2,\ldots)\in\J^{(e)}$.
If $\tilde\nu$ is a non-standard generalized composition, then let
$\tilde\alpha=\tilde\pi\tilde\rho$ be the factorization such that
$ d(\tilde\pi)-\ell(\tilde\pi) \ge e$, and  let
$\tilde\alpha=\tilde\nu\,0\,0\cdots=\tilde\gamma\, 0\, a\,\beta\,0\,0\cdots$ be the
factorization used in~(\ref{inductG}). If $\tilde\pi$ is an initial factor of
$\tilde\gamma$, then it is clearly an initial  factor of both ${\tilde\gamma\, 
a\,\beta\,0\,0\cdots}$ and
${\tilde\gamma\, (a-1)\,\beta\,0\,0\cdots}$ and they both reach level $e$. By
the induction hypothesis, both 
$G_{\tilde\gamma\,  a\,\beta\,0\cdots}$ and  $G_{\tilde\gamma\, (a-1)\,\beta\,0\cdots}$
are in $\J^{(e)}$ and in turn
$G_{\tilde\alpha}\in\J^{(e)}$. If we now assume that $\tilde\pi=\tilde\gamma\,0$, then
  $$ d(\tilde\gamma)-\ell(\tilde\gamma)=d(\tilde\pi)-\big(\ell(\tilde\pi)-1\big)\ge e+1>e$$
and again the induction hypothesis can be applied to~(\ref{inductG}) to show that
$G_{\tilde\alpha}\in\J^{(e)}$.

 We are left to check the case where
  $$\tilde\pi=\tilde\gamma\,0\,a\,\tilde\mu.$$
For the first term in~(\ref{inductG}), $\tilde\gamma\,a\,\tilde\mu$ is an initial
factor and
  $$d(\tilde\mu\,a\,\tilde\gamma)-\ell
(\tilde\gamma\,a\,\tilde\mu)=d(\tilde\pi)-\big(\ell(\tilde\pi)-1\big)
    \ge e+1>e .$$
The induction hypothesis gives that $G_{\tilde\gamma\, 
a\,\beta\,0\cdots}\in\J^{(e)}$. 
For the second term indexed by $\tilde\gamma\, (a-1)\,\beta\,0\,0\cdots$
we have 
  $$d(\tilde\gamma(a-1)\tilde\mu)-\ell(\tilde\gamma(a-1)\tilde\mu)=
   \big(d(\tilde\pi)-1\big)-\big(\ell(\tilde\pi)-1\big)\ge e.$$ 
 Again the induction hypothesis gives us that $G_{\tilde\gamma\,
(a-1)\,\beta\,0\cdots}\in\J^{(e)}$, concluding the proof that $G_{\tilde\alpha}\in\J^{(e)}$.

For the second statement of the lemma let ${\tilde\gamma\,(a-1)\,\beta\,0\,0\cdots
}=\tilde\pi\tilde\rho$ be a factorization such that $d(\tilde\pi)-\ell(\tilde\pi)\ge e$.
If $\tilde\pi$ is an initial factor of
$\tilde\gamma$ then it is clear that $\tilde\alpha$ reaches level $e$. On the
other hand if $\tilde\pi=\tilde\gamma\,(a-1)\,\tilde\mu$, then we have
  $$d(\tilde\gamma\,0\,a\,\tilde\mu)-\ell(\tilde\gamma\,0\,a\,\tilde\mu)=
    \big(d\big(\tilde\gamma\,(a-1)\,\tilde\mu\big)+1\big)-
   \big(\ell\big(\tilde\gamma\,(a-1)\,\tilde\mu\big)+1\big)
   \ge e.$$
 Thus $\tilde\alpha$ reaches level $e$ which concludes our proof.
\endproof

In light of the previous lemma, let 
${\mathcal G}^{(e)}$  
denote the set of all
generalized infinite compositions
$\tilde\alpha$ reaching level $e$, that satisfy~(\ref{econdition}). We remark that
the set
$\{ G_{\tilde\alpha}\}_{\tilde\alpha\in{\mathcal G}^{(e)}}$  constructed above is
contained in $\J^{(e)}$ and contains $\{F_\alpha(x_1,x_2,\ldots)\,:\,\exists\,\alpha=
\pi\rho,\ d(\pi)-\ell(\pi)\ge e\}$. Hence we have

  \begin{lem} \label{idealG}
   $$\J^{(e)} = \langle G_{\tilde\alpha}  \rangle_{\tilde\alpha\in{\mathcal
G}^{(e)}}\ .$$
\endproof
  \end{lem} 

Our next task is to characterize the leading monomial of each function $G_{\tilde\alpha}$.
Before this we need to specify which monomial order we use.
 Let
$X^{\tilde\alpha}=x_1^{\tilde\alpha_1}x_2^{\tilde\alpha_2}\cdots$ and
$X^{\tilde\beta}=x_1^{\tilde\beta_1}x_2^{\tilde\beta_2}\cdots$ be any two
monomials where $\tilde\alpha$ and $\tilde\beta$ are two generalized infinite
compositions. We say that
$X^{\tilde\alpha}\le_{lex} X^{\tilde\beta}$ if and only if
$d(\tilde\alpha)>d(\tilde\beta)$, or
$d(\tilde\alpha)=d(\tilde\beta)$ and the leftmost non-zero entry in
$[{\tilde\beta_1}-{\tilde\alpha_1},{\tilde\beta_2}-{\tilde\alpha_2},
\ldots,]$ is positive. The order $\le_{lex}$ is a classical monomial
order in the sense that it is a total order and if $X^{\tilde\alpha}\le_{lex}
X^{\tilde\beta}$, then $X^{\tilde\alpha}X^{\tilde\gamma}=
X^{\tilde\alpha+\tilde\gamma}\le_{lex} X^{\tilde\beta+\tilde\gamma}=
X^{\tilde\beta}X^{\tilde\gamma}$. Here the sum of generalized  compositions is
componentwise.

For any formal power series $P=P(x_1,x_2,\ldots)\in\Q[[x_1,x_2,\ldots]]$  we let
$LM(P)$ denote the leading monomial of
$P$. That is $LM(P)$ is the monomial of $P$ with non-zero coefficient  of smallest
degree and largest in lexicographic order. In other words, the leading monomial for the
order
$\le _{lex}$. We let
$LC(P)$ denote the coefficient of $LM(P)$ in $P$. Remark  that for any two functions
$P$ and  $Q$, we have
  $$LM (P Q)= LM(P) LM(Q).$$

We need the following result which is the extension  for the $G$-functions of 
Relations~(\ref{firstFrel}) and~(\ref{secondFrel}) for the $F$-functions.
  \begin{lem} \label{Lessvar} 
      Let $\tilde\alpha=b\,\tilde\rho$ be any  generalized infinite composition and
$b\ge 0$.
   \begin{enumerate}
     \item If $b=0$, then $G_{\tilde\alpha}$
  is a functions with no variable $x_1$. More precisely, we have
  \begin{equation}\label{firstGrel} G_{0\,\tilde\mu}(x_1,x_2,\ldots)=
G_{\tilde\rho}(x_2,x_3,\ldots).\end{equation}
  \item If $b>0$, then
  \begin{equation}\label{secondGrel} G_{\tilde\alpha}= x_1 G_{(b-1)\,\tilde\rho}\  +\
     M_{\tilde\alpha}(x_2,x_3,\ldots),\end{equation}
  where $M_{\tilde\alpha}\in\J^{(e)}(x_2,x_3,\ldots)$  whenever
$G_{\tilde\alpha}\in\J^{(e)}$
  \end{enumerate}\end{lem} 

\proof  Remark that in~(\ref{firstGrel}),  the generalized composition $\tilde\rho$ is
also infinite and\break
$G_{\tilde\rho}(x_2,x_3,\ldots)$ is the function
$G_{\tilde\rho}(x_1,x_2,\ldots)$ in which the variable $x_i$ is replaced by
$x_{i+1}$.  Similarly, $\J^{(e)}(x_2,x_3,\ldots)$  is the ideal $\J^{(e)}$ where each 
variable
$x_i$ is replaced by $x_{i+1}$.

To show the two relations, we let $\tilde\alpha= \tilde\nu\,0\,0\cdots$ and proceed by
induction on
$\ell(\tilde\nu)$. If $\ell(\tilde\nu)=0$ then we have
$G_{0\,0\cdots}=1$ and~(\ref{firstGrel}) is valid.
In the following we extensively use the two relations~(\ref{firstFrel})
and~(\ref{secondFrel}).  If
$\tilde\alpha=b\,\rho\,0\,0\cdots$, then Equation~(\ref{baseG}) gives us
$G_{\tilde\alpha}=F_{b\,\rho}(x_1,x_2,\ldots)$. If $b>1$, then we get
  $$\begin{array}{rcl}
   G_{b\,\rho\,0\cdots}&=&F_{b\,\rho}(x_1,x_2,\ldots) \ =\ x_1
F_{(b-1)\,\rho}(x_1,x_2,\ldots)\  +\ F_{b\,\rho}(x_2,x_3,\ldots)\\
   &=&x_1 G_{(b-1)\,\rho\,0\cdots}\ +\
        M_{b\,\rho\,0\cdots}(x_2,x_3,\ldots),
   \end{array}
  $$
and (\ref{secondGrel}) follows for this case with
$M_{\tilde\alpha}(x_2,x_3,\ldots)=F_{b\,\rho}(x_2,x_3,\ldots)$.

 For $b=1$, we first need to understand (\ref{firstGrel}) in the case
$G_{0\,\rho\,0\cdots}$.  For this assume that  $\rho=a\beta$. If $a>1$, then the
Definitions~(\ref{baseG}) and~(\ref{inductG})  give
  $$\begin{array}{rcl}
   G_{0\,a\,\beta\,0\cdots}&=&G_{a\,\beta\,0\cdots}-x_1 G_{(a-1)\,\beta\,0\cdots} \ =\
F_{a\,\beta}(x_1,x_2,\ldots)-x_1 F_{(a-1)\,\beta}(x_1,x_2,\ldots)\\
  & =&
    F_{a\,\beta}(x_2,x_3,\ldots)=G_{a\,\beta\, 0\cdots}(x_2,x_3,\ldots)\\
   \end{array}
  $$
If $a=1$, then we use the induction hypothesis on 
$\ell(0\,\beta)=\ell(0\,1\,\beta)-1<\ell(\tilde\nu)$ to get
  $$\begin{array}{rcl}
   G_{0\,1\,\beta\,0\cdots}&=&G_{1\,\beta\,0\cdots}-x_1 G_{0\,\beta\,0\cdots} 
   \ =\ F_{1\,\beta}(x_1,x_2,\ldots)-x_1 F_{\beta}(x_2,x_3,\ldots)\\
   & =&
    F_{1\,\beta}(x_2,x_3,\ldots)=G_{1\,\beta \, 0\cdots}(x_2,x_3,\ldots)\\
   \end{array}
  $$
Now we can go back to (\ref{secondGrel}) in the  case of
$\tilde\alpha=1\,\rho\,0\cdots$:
  $$\begin{array}{rcl}
   G_{1\,\rho\,0\cdots}&=&F_{1\,\rho}(x_1, x_2,\ldots)\ =\ x_1
F_{\rho}(x_2,x_3,\ldots)\  +\ F_{1\,\rho}(x_2,x_3,\ldots)\\
   &=&x_1 G_{0\,\rho\,0\cdots}\ +\
        M_{1\,\rho\,0\cdots}(x_2,x_3,\ldots).
   \end{array}
  $$
We remark that in both cases, for $\tilde\alpha=\nu 0\cdots$, we have 
  \begin{equation}\label{baseM}
     M_{\nu\,0\cdots}(x_2,x_3,\ldots)=F_{\nu} (x_2,x_3,\ldots)\in\J^{(e)}(x_2,x_3\ldots)
  \end{equation}
whenever $G_{\nu\,0\cdots}\in\J^{(e)}(x_2,x_3\ldots)$

We then consider when $\tilde\alpha=\tilde\gamma\,  0\, a\,\beta\,0\,0\cdots$. This is
the factorization needed to use (\ref{inductG}) with
$k=\ell(\tilde\gamma)+1$. If $\tilde\gamma$ is empty,  then we have $b=0$ and we are in
the case considered above. Assume  that $\tilde\gamma=b\,\tilde\mu$. If $b=0$, then
applying the induction hypothesis we have
  $$\begin{array}{rcl}
    G_{0\,\tilde\mu\,0\,  a\,\beta\,0\cdots}(x_1,x_2,\ldots)
     &=& G_{0\,\tilde\mu\,  a\,\beta\,0\cdots}(x_1,x_2,\ldots) - x_k
G_{0\,\tilde\mu\,(a-1)\,\beta\,0\cdots}(x_1,x_2,\ldots)\cr
     &=&G_{\tilde\mu\,  a\,\beta\,0\cdots}(x_2,x_3,\ldots) - x_{(k-1)+1}
             G_{\tilde\mu\,(a-1)\,\beta\,0\cdots}(x_2,x_3,\ldots)\\
     &=&G_{\tilde\mu\,0\,  a\,\beta\,0\cdots}(x_2,x_3,\ldots).\end{array}$$
Here remark that even though $\ell(\tilde\mu)+1=k-1$,  we have to replace $x_{k-1}$ by
$x_{(k-1)+1}=x_k$ in the defining recurrence for $G_{\tilde\mu\,0\, 
a\,\beta\,0\cdots}(x_2,x_3,\ldots)$. If
$b>0$, the induction hypothesis now gives
   $$\begin{array}{rcl}
    G_{b\,\tilde\mu\,0\,  a\,\beta\,0\cdots}
     &=& G_{b\,\tilde\mu\,  a\,\beta\,0\cdots}  - x_k
G_{b\,\tilde\mu\,(a-1)\,\beta\,0\cdots}\cr
     &=&x_1 G_{(b-1)\,\tilde\mu\,  a\,\beta\,0\cdots} + M_{b\,\tilde\mu\, 
       a\,\beta\,0\cdots}(x_2,x_3,\ldots)\\
      &&\  -\  x_k \big(x_1 G_{(b-1)\,\tilde\mu\,(a-1)\,\beta\,0\cdots}+
            M_{b\,\tilde\mu\,(a-1)\,\beta\,0\cdots}(x_2,x_3,\ldots)\big)\cr
     &=&x_1 \big( G_{(b-1)\,\tilde\mu\,  a\,\beta\,0\cdots}-x_k
          G_{(b-1)\,\tilde\mu\,(a-1)\,\beta\,0\cdots}\big)\\
    &&\  +\   M_{b\,\tilde\mu\, a\,\beta\,0\cdots}(x_2,x_3,\ldots)
       - x_k  M_{b\,\tilde\mu\,(a-1)\,\beta\,0\cdots}(x_2,x_3,\ldots)\cr
     &=&x_1 G_{(b-1)\,\tilde\mu\,0\,  a\,\beta\,0\cdots } + M_{b\,\tilde\mu\,0\, 
        a\,\beta\,0\cdots}(x_2,x_3,\ldots),
      \cr
   \end{array}$$
where the function 
  \begin{equation}\label{inductM}
     M_{b\,\tilde\mu\,0\,  a\,\beta\,0\cdots}(x_2,x_3,\ldots) =\big(M_{b\,\tilde\mu\,
a\,\beta\,0\cdots}
       - x_k  M_{b\,\tilde\mu\,(a-1)\,\beta\,0\cdots}\big)(x_2,x_3,\ldots)
  \end{equation}
contains no variable $x_1$. Using the argument in  Lemma~\ref{InJ}, if
$G_{b\,\tilde\mu\,0\,  a\,\beta\,0\cdots}\in\J^{(e)}$, then both $G_{b\,\tilde\mu\,
a\,\beta\,0\cdots}$ and
$G_{b\,\tilde\mu\,(a-1)\,\beta\,0\cdots}$ are in
$\J^{(e)}$. The induction hypothesis gives us that $M_{b\,\tilde\mu\,0\, 
a\,\beta\,0\cdots}(x_2,x_3,\ldots)\in\J^{(e)}(x_2,x_3,\ldots)$  and this completes the
proof of the lemma.
\endproof

  \begin{cor} \label{LMG} 
   Let $\tilde\alpha$ be any generalized infinite composition. We
have
   $$LM(G_{\tilde\alpha})= X^{\tilde\alpha}.$$ 
  \end{cor} 
\proof Let $\tilde\alpha=\tilde\nu\,0\,0\cdots$. We proceed  by induction on
$\ell(\tilde\nu)$ and the degree
$d=d(\tilde\alpha)=\sum\tilde\alpha_i$. If $\ell(\tilde\nu)=0$  we have
$G_{0\,0\,\cdots}=1=X^{0\,0\cdots}$. If
$\ell(\tilde\nu)\ge1$, then let
$\tilde\alpha=b\,\tilde\rho $ as in Lemma~\ref{Lessvar}.  If $b=0$, then the induction
hypothesis on
$\ell(\tilde\nu)$ gives
  $$LM(G_{0\,\tilde\rho })= LM\big(G_{\tilde\rho}(x_2,x_3,\ldots)\big)=
  x_1^0x_2^{\tilde\rho_1}x_3^{\tilde\rho_2}\cdots =X^{\tilde\alpha}.$$ 
Now if $b>0$ we
use the second part of Lemma~\ref{Lessvar} and the induction hypothesis on $d$, and get
  $$LM\big(x_1 G_{ (b-1)\,\tilde\mu} + M_{\tilde\alpha}(x_2,x_3,\ldots)\big)=x_1 LM(G_{
      (b-1)\,\tilde\mu})= X^{\tilde\alpha}.
  $$
\endproof

\begin{rem}\label{remarque} \rm From the above corollary,  by triangularity, it is
clear that the set
$\{G_{\tilde\alpha}\}$ for all $\tilde\alpha$ forms a  Hilbert basis of
$\Q[[x_1,x_2\ldots]]$. We will see in Section~5 that in fact
$\{G_{\tilde\alpha}\}_{\tilde\alpha\in{\mathcal G}^{(e)}}$  forms a Hilbert basis of
$\overline{\J^{(e)}}$.
\end{rem}

%%%%%%%%%%%%%%%%%%%%%%%%%%%%
\section{It is at most Catalan.}

Let ${\mathcal Q}^{(e)}=\{ G_{\tilde\alpha} \}_{\tilde\alpha\in{\mathcal G}^{(e)}}$ be
the generating set of $\J^{(e)}$ constructed  in Lemma~\ref{idealG}. In this section we
show that after reduction, at most the monomials corresponding to Catalan paths form a
Hilbert basis of $R =\Q[[x_1,x_2,\ldots]]\big/\lower2pt\hbox{$\overline{\J}$}$. For
this we reduce every other monomial to these. In fact for $R^{(e)}=\Q[[x_1,x_2,\ldots]]
\big/ \lower3pt\hbox{$\overline{\J^{(e)}}$}$,  we show that at most the monomials
corresponding to paths above the line $y=x-e$ form a Hilbert basis of
$R^{(e)}$, for all
$e\ge 0$. We conclude this section with the corresponding result for finitely many
variables, $R_n^{(e)}=\Q[x_1,x_2,\ldots,x_n]
\big/ \lower3pt\hbox{${\J_n^{(e)}}$}$, which  is a generalization of
Theorem~\ref{Finite}.

Given any generalized infinite composition $\tilde\alpha$  we associate a
unique path in the plane with steps going north or east. More
precisely, for
$\tilde\alpha=(\tilde\alpha_1,\tilde\alpha_2,\tilde\alpha_3,\ldots)$,  we construct the
path that starts at
$(0,0)$, then moves $\tilde\alpha_1$ steps east to $(\tilde\alpha_1,0)$;  then one step
north to
$(\tilde\alpha_1,1)$, and then $\tilde\alpha_2$ steps east to 
$(\tilde\alpha_1+\tilde\alpha_2,1)$; then one step north to
$(\tilde\alpha_1+\tilde\alpha_2,2)$, and then $\tilde\alpha_3$ steps east to
$(\tilde\alpha_1+\tilde\alpha_2+\tilde\alpha_3,2)$; and so on.
For example for $\tilde\alpha=0\,0\,2\,1\,0\,0\,3\,0\,0\cdots$ 
we have the path
    $$\Dyck{&&&&&&&\G\cr
            &&&&\Bas{$x_7$}&\Bas{$x_7$}&\Bas{$x_7$}&\G\cr
            &&&&\G\cr
            &&&&\G\cr
            &&&\Bas{$x_4$}&\G\cr
            &\Bas{$x_3$}&\Bas{$x_3$}&\G\cr
            &\G\cr
            &\G\cr
     }$$
For every east step at hight $i-1$ we associate a variable $x_i$.  The product of all
the variables associate to a path encoded by $\tilde\alpha$ is denoted
$X^{\tilde\alpha}$. We now remark that for any factorization
$\tilde\alpha=\tilde\pi\tilde\rho$, the rightmost coordinate of the path at hight
$\ell(\tilde\pi)-1$ is $\big(d(\tilde\pi),\ell(\tilde\pi)-1\big)$.
 \begin{definition} For an integer $e\ge 0$, we say that a generalized composition
$\tilde\alpha$  is of type
$e$-Catalan if 
 its associated path remains above the line $y=x-e$. That is,  every coordinate
$(x_i,y_i)$ of the path 
 is such that $x_i-y_i\le e$.
  \end{definition}
  \begin{lem} \label{upper} The monomials of $\Q[[x_1,x_2,\ldots]]$  corresponding to
paths remaining above the line $y=x-e$ contains a Hilbert basis of the quotient
$R^{(e)}$.
  \end{lem} 

\proof Let
$X^{\tilde\alpha}$ be any monomial of degree $d$. If the path  corresponding to
$\tilde\alpha$ goes under the line
$y=x-e$, then let $\tilde\alpha=\tilde\pi\tilde\rho$ be any factorization  such that
the coordinate
$\big(d(\tilde\pi),\ell(\tilde\pi)-1\big)$ is under the line $y=x-e$. That is
  $$d(\tilde\pi)-\ell(\tilde\pi)\ge e.
   $$
From Lemma~\ref{InJ} we conclude that the
function
$G_{\tilde\alpha}$ with leading monomial
$X^{\tilde\alpha}$ is in $\J^{(e)}$. This monomial can thus be replaced  by monomials
of degree $d$ but strictly smaller with respect to
$<_{lex}$. 
Repeating this step (possibly countably many times) with the next largest monomial 
going under the line
$y=x-e$, any monomial
$X^{\tilde\alpha}$ can be reduced modulo the ideal
$\overline{\J^{(e)}}$ to a series containing only monomials
$X^{\tilde\beta}$ where $\tilde\beta$ is of type $e$-Catalan.
\endproof

We are now in position to generalize Theorem~\ref{Finite} and prove it.
For this, note that the quasi-symmetric polynomials in $n$ variables are defined by
setting $0=x_{n+1}=x_{n+2}=\cdots$ in the quasi-symmetric functions. That is
  $$ F_\alpha(x_1,x_2,\ldots,x_n)=F_\alpha(x_1,x_2,\ldots,x_n,0,0,\ldots).
  $$
We then define 
 $$J_n^{(e)}=\langle F_\alpha(x_1,x_2,\ldots,x_n)\,:\,\alpha \hbox{ \sl
    reaches level }e\rangle$$ 
and $R_n^{(e)}=\Q[x_1,x_2,\ldots,x_n]\big/J_n^{(e)}$.
Similarly we set
  $$
G_{\tilde\alpha}(x_1,x_2,\ldots,x_n)=G_{\tilde\alpha}(x_1,x_2,\ldots,x_n,0,0,\ldots).
  $$
It is clear that Lemma~\ref{InJ} holds for $J_n^{(e)}$ in the same way. More over if
$\tilde\alpha=\tilde\nu\,0\,0\cdots$ for $\ell(\tilde\nu)=n$ then 
 $$LM(G_{\tilde\alpha}(x_1,x_2,\ldots,x_n))=
   LM(G_{\tilde\alpha}(x_1,x_2,\ldots,x_n,0,0,\ldots))=x_1^{\tilde\nu_1}
   x_2^{\tilde\nu_2} \cdots x_n^{\tilde\nu_n}.
 $$
Let $C_n^{(e)}$ denotes the number
of generalized compositions
$\tilde\alpha$ of type $e$-Catalan such that $\tilde\alpha=\tilde\nu\,0\,0\cdots$ and
$\ell(\tilde\nu)=n$.
These are in bijection with the paths from  $(0,0)$ to $(n+e,n)$ that remain above
the line $y=x-e$. Indeed, if we have a path of type $e$-Catalan, it suffices to add a
horizontal line from $(\tilde\nu_1+\tilde\nu_2+\cdots+\tilde\nu_n,n)$ to $(n+e,n)$.  
When $e=0$, we have $C_n^{(0)}=C_n$ the $n$th Catalan number. This enumerates the
classical Dyck path from $(0,0)$ to
$(n,n)$ remaining above the line
$y=x$. See
\cite{Stanley_EC1} for an extensive account on Catalan numbers. We have the following
generalization to our Theorem~\ref{Finite}.

  \begin{cor} \label{Fupper} $\dim(R_n^{(e)})\le C_n^{(e)}$.
  \end{cor} 
\proof We use the same argument as in Lemma~\ref{upper}. For this we use the fact  that
for any monomial
$x_1^{\tilde\nu_1}x_2^{\tilde\nu_2}\cdots x_n^{\tilde\nu_n}$ in
$\Q[x_1,x_2,\ldots,x_n]$, if $\tilde\nu$ reaches level $e$, then
$G_{\tilde\nu\,0\cdots}(x_1,x_2,\ldots,x_n)\in\J_n^{(e)}$. Hence a basis of $R_n^{(e)}$
is contained in the monomials corresponding to paths of type $e$-Catalan and our result
follows.
\endproof

Again, we expect the equality to hold in Corollary~\ref{Fupper}, and we will address
this question in \cite{deux}.

%%%%%%%%%%%%%%%%%%%%%%%%%%%%
\section{It is a Hilbert basis}
 
In the previous Section, the generating set ${\mathcal Q}^{(e)}=\{
G_{\tilde\alpha}\}_{\tilde\alpha\in{\mathcal G}^{(e)}}$ of $\J^{(e)}$  is very useful
to reduce every monomial to $e$-Catalan type generalized compositions. It is in fact a
Hilbert basis for the given ideal. We use here ideas of Gr\"obner basis theory. This
is crucial to complete the proof of Theorem~\ref{main}.  Let us recall a few basic
facts about Gr\"obner bases, see
\cite{grob,cox} for more details. 

To show that a set $S$ is a Gr\"obner basis
it is enough to show that all polynomial syzygies of that  set are reducible in $S$.
The polynomial syzygy of
$P$ and $Q$ is defined by
  \begin{equation} \label{syzy} 
       S(P,Q)=LC(Q) M_1 P - LC(P) M_2 Q
  \end{equation} 
where $\hbox{lcm}\big(LM(P),LM(Q)\big)= M_1\cdot LM(P) = M_2\cdot LM(Q)$.  This shows
that the given set contains all the generators of the leading monomials of the ideal.

To help us we use the classic Buchberger's lemma \cite{grob,cox}:
  \begin{lem}\label{Buch}
    Given $P,Q\in S$. If there is an $R\in S$ such that $LM(R)$
    divides\break $\hbox{lcm}\big(LM(P),LM(Q)\big)$, and if both  $S(R,Q)$ and $S(R,P)$
are reducible in
    $S$, then $S(P,Q)$ is reducible in $S$. \endproof
  \end{lem}
This result is easily adapted to our context. We first remark that  our sets ${\mathcal
Q}^{(e)}$ are {\sl lattice}. That is 
  \begin{lem}\label{latt}
  $$ G_{\tilde\alpha}\in{\mathcal Q}^{(e)} \qquad\implies\qquad 
G_{\tilde\rho}\in{\mathcal Q}^{(e)}$$ for all ${\tilde\alpha}\le {\tilde\rho}$
componentwise.
  \end{lem}
\proof
If $\tilde\alpha=\tilde\pi\tilde\nu$ satisfies  $d(\tilde\pi)-\ell(\tilde\pi)\ge e$,
then let
$r=\ell(\tilde\pi)$ and consider $\tilde\rho=\tilde\gamma\tilde\mu$  where
$\ell(\tilde\gamma)=r$. Since $d(\tilde\pi)\le d(\tilde\gamma)$, we have
  $$d(\tilde\gamma)-\ell(\tilde\gamma)\ge d(\tilde\pi)-\ell(\tilde\pi) \ge e.
 $$
By Lemma~\ref{InJ}, $G_{\tilde\rho}\in{\mathcal Q}^{(e)}$.
\endproof

We can now adapt the proof (see \cite{grob,cox}) of Lemma~\ref{Buch}  to our
situation. For any pair
${\tilde\alpha},{\tilde\pi}\in{\mathcal G}^{(e)}$, we define 
$S(G_{\tilde\alpha},G_{\tilde\pi})$ as in (\ref{syzy}) with our  definition of $LM$ and
$LC$. We show in this section that any such $S(G_{\tilde\alpha},G_{\tilde\pi})$  is
reducible in
${\mathcal Q}^{(e)}$. Let
$\tilde\rho\in{\mathcal G}^{(e)}$ be the unique element such that 
    $$X^{\tilde\rho}=
      \hbox{lcm}\big(X^{\tilde\alpha},
          X^{\tilde\pi}\big)= M_1 X^{\tilde\alpha}= M_2X^{\tilde\alpha}.$$ 
We have
  \begin{equation}\label{Buch1}
  \begin{array}{rcl}
    S(G_{\tilde\alpha},G_{\tilde\pi}) &=&  M_1 G_{\tilde\alpha} - M_2 G_{\tilde\pi}\cr
   &=& M_1 G_{\tilde\alpha}- G_{\tilde\rho} +G_{\tilde\rho} - M_2 G_{\tilde\pi}\cr
   &=& S(G_{\tilde\alpha},G_{\tilde\rho}) + S(G_{\tilde\rho},G_{\tilde\pi}). \cr
    \end{array}
  \end{equation}
If both $S(G_{\tilde\alpha},G_{\tilde\rho})$ and $S(G_{\tilde\rho},G_{\tilde\pi})$  
are  reducible in
${\mathcal Q}^{(e)}$, then so is $ S(G_{\tilde\alpha},G_{\tilde\pi})$.  It is thus
sufficient to show that all
$S(G_{\tilde\alpha},G_{\tilde\rho})$ are reducible in ${\mathcal Q}^{(e)}$ for
$\tilde\alpha\le\tilde\rho$ componentwise.

We can reduce further our problem as follows. Assume that   $\tilde\alpha$ and
$\tilde\rho$ in ${\mathcal G}^{(e)}$ are generalized compositions of $d_1$  and $d_2$
respectively. If
$\tilde\alpha\le\tilde\rho$, then $d_1\le d_1$. If $d_2-d_1>1$,  we can select a
generalized composition 
$\tilde\alpha\le\tilde\pi\le\tilde\rho$ and use (\ref{Buch1})  again. We can thus
assume that
$d_2-d_1=1$. That is the two generalized compositions differ  on one part only and by
one unit.

  \begin{lem}\label{Gbasis}
    The set ${\mathcal Q}^{(e)}$ is a Gr\"obner basis of $\J^{(e)}$.
  \end{lem}
\proof
From the discussion above it is sufficient to show that all  the expressions of
${\mathcal Q}^{(e)}$ of the form 
  \begin{equation} \label{syzy1}
       S(G_{\tilde\gamma\,a\,\tilde\beta},G_{\tilde\gamma\,(a-1)\,\tilde\beta})=
       G_{\tilde\gamma\,a\,\tilde\beta}- x_k G_{\tilde\gamma\,(a-1)\,\tilde\beta}
  \end{equation} 
where $k=\ell(\tilde\gamma)+1$,  are reducible in ${\mathcal Q}^{(e)}$. Let us denote by
$m_{\tilde\gamma\,a\,\tilde\beta}(x_1,x_2,\ldots)$ the leading monomial of
$S(G_{\tilde\gamma\,a\,\tilde\beta},G_{\tilde\gamma\,(a-1)\,\tilde\beta})$.

Let $\tilde\beta=\tilde\nu\,0\,0\cdots$. We set up an induction on $\ell(\tilde\nu)$.
Assume first that $\tilde\nu=\nu$ is a (possibly empty) standard  composition. The
second part of Lemma~\ref{InJ} and the recursive definition~(\ref{inductG}) give
  \begin{equation} \label{syzy2}
     S(G_{\tilde\gamma\,a\,\nu\,0\cdots},
      G_{\tilde\gamma\,(a-1)\,\nu\,0\cdots})=
     G_{\tilde\gamma\,0\,a\,\nu\,0\cdots}\quad \in\ {\mathcal Q}^{(e)}.
  \end{equation}  

If $\tilde\nu$ is not standard, then let
$\tilde\beta=\tilde\pi\,0\,b\,\mu\,0\,0\cdots $ for $b>0$  and $\mu$  a (possibly
empty) standard composition. Let $\ell=\ell(\tilde\gamma\,a\,\tilde\pi)+1$.
Using~(\ref{inductG}), we have
  \begin{equation} \label{syzy2.5}
    \begin{array}{rcl}
     S(G_{\tilde\gamma\, a\,\tilde\pi\,0\,b\mu0\cdots},
       G_{\tilde\gamma(a-1)\tilde\pi\,0\,b\mu0\cdots})
    &=&
       G_{\tilde\gamma\,a\,\tilde\pi\,0\,b\,\mu\,0\cdots}- x_k
       G_{\tilde\gamma\,(a-1)\,\tilde\pi\,0\,b\,\mu\,0\cdots}.\\
    &=& G_{\tilde\gamma\,a\,\tilde\pi\,b\,\mu\,0\cdots}-x_\ell
       G_{\tilde\gamma\,a\,\tilde\pi\,(b-1)\,\mu\,0\cdots}\\
    &&-
       x_k
       \big(G_{\tilde\gamma(a-1)\tilde\pi\,b\,\mu0\cdots}- x_\ell
         G_{\tilde\gamma(a-1)\tilde\pi(b-1)\mu0\cdots}\big)\\
    &=& S(G_{\tilde\gamma\,a\,\tilde\pi\,b\,\mu\,0\cdots},
      G_{\tilde\gamma\,(a-1)\,\tilde\pi\,b\,\mu\,0\cdots})\\
    &&- 
    x_\ell\, S(G_{\tilde\gamma\,a\,\tilde\pi(b-1)\mu\,0\cdots},
     G_{\tilde\gamma(a-1)\tilde\pi(b-1)\mu\,0\cdots}).
     \end{array}
  \end{equation}  
We observe
that
  \begin{equation} \label{syzym}
     m_{\tilde\gamma\, a\,\tilde\pi\,0\,b\mu0\cdots}(x_1,x_2,\ldots)= x_\ell\, 
     m_{\tilde\gamma\,a\,\tilde\pi(b-1)\mu\,0\cdots}(x_1,x_2,\ldots)
  \end{equation}  
since $m_{\tilde\gamma\,a\,\tilde\pi\,b\,\mu\,0\cdots}(x_1,x_2,\ldots)$  is distinct
and lexicographically smaller than\break
$x_\ell\,  m_{\tilde\gamma\,a\,\tilde\pi(b-1)\mu\,0\cdots}(x_1,x_2,\ldots)$.  Thus the
left hand side of (\ref{syzy2.5}) is resolved by two expressions such that
$\ell(\tilde\pi\,b\,\mu)=\ell(\tilde\pi(b-1)\mu)<\ell(\tilde\nu)$  and by the induction
hypothesis the right hand side of (\ref{syzy2.5}) is reducible in ${\mathcal Q}^{(e)}$.
This completes the proof of the lemma.
\endproof

  \begin{cor}\label{Hbasis}
    The set ${\mathcal Q}^{(e)}$ is a Hilbert basis of $\overline{\J^{(e)}}$.
  \end{cor}
\proof
Given an element $P\in\overline{\J^{(e)}}$ let $X^{\tilde\beta}=LM(P)$. Since 
${\mathcal Q}^{(e)}$ is a Gr\"obner basis it contains an element $G_{\tilde\alpha}$
such that
$\tilde\alpha \le \tilde\beta$ componentwise. Lemma~\ref{latt} gives us that
$G_{\tilde\beta}\in{\mathcal Q}^{(e)}$. The element
$P-LC(P)\cdot G_{\tilde\alpha}\in\overline{\J^{(e)}}$ is such that $LM(P-LC(P)\cdot
G_{\tilde\alpha})<_{lex} X^{\tilde\beta}$. If we repeat this process (possibly
countably many times) we can express $P$ as a series in the elements of
${\mathcal Q}^{(e)}$.
\endproof

We are now in position to conclude our investigation and  prove the general version of
our Theorem~\ref{main}. 

  \begin{cor}\label{Maine}
    the monomial Hilbert basis of $R^{(e)}$ is given by the monomials of $\Q[[x_1,
    x_2,x_3,\ldots]]$ corresponding to the paths of type $e$-Catalan.
  \end{cor}
\proof
 As noted in Remark~\ref{remarque}, the set $\{G_{\tilde\alpha}\}$  forms a Hilbert
basis of
$\Q[[x_1,x_2,x_3,\ldots]]$, and Corollary~\ref{Hbasis} gives that  the set ${\mathcal
Q}^{(e)}$ is a Hilbert basis of $\overline{\J^{(e)}}$. Thus a Hilbert basis of the
quotient
$R^{(e)}=\Q[[x_1,x_2,x_3,\ldots]]\big/\lower3pt\hbox{$\overline{\J^{(e)}}$}$ is given by
the set 
 $$\{G_{\tilde\alpha}\}\setminus {\mathcal Q}^{(e)}=\big \{G_{\tilde\alpha}\ \big|\
{\tilde\alpha}
    \hbox{\ \sl corresponds to a path of type $e$-Catalan}
   \big\}.$$
The result follows by triangularity.\endproof

  \begin{rem}\rm
    To show the equality in Equation~(\ref{dSPACE}), it appears that the set
$$\{G_{\tilde\alpha}\,|\,{\ell(\tilde\alpha)=n},\ \hbox{\sl and $\tilde\alpha$
 reaches level $e$}\}$$ 
forms a linear basis of $J_n^{(e)}$. Unfortunately the argument
of Section~5 is not sufficient to show this with finitely many variables and it requires
more work. This is the object of our collaboration in \cite{deux}.
  \end{rem}
%

%%%%%%%%%%%%%%%%%%%%%%%%%%%%
\section*{Acknowledgments}
We are indebted to Fran\c{c}ois Bergeron and Adriano  Garsia for extremely
stimulating conversations and valuable suggestions. We  also thank Mike Zabrocki and
Geanina Tudose for ideas and comments.

%%%%%%%%%%%%%%%%%%%%%%%%%%%%

\end{document}